\documentclass[11pt]{article}
\usepackage{latexsym}
\usepackage{amssymb,amsmath}
\usepackage{graphicx}
\RequirePackage{tikz}
\usepackage{epsfig}
\usepackage{amsfonts}
\newtheorem{theorem}{Theorem}[section]

\newtheorem{proposition}[theorem]{Proposition}

\newcommand{\proof}{\noindent{\bf Proof.\ }}
\newcommand{\qed}{\hfill $\square$ \bigskip}

\newcommand{\edim}{{\rm edim}}

\begin{document}

\title{Edge metric dimensions via hierarchical product and integer linear programming}

\author{Sandi Klav\v zar $^{a,b,c}$ \and Mostafa Tavakoli $^{d,}$\footnote{Corresponding author}} 

\date{}

\maketitle

\vspace*{-10mm}
\begin{center}
$^a$ Faculty of Mathematics and Physics, University of Ljubljana, Slovenia\\
{\tt sandi.klavzar@fmf.uni-lj.si}

\medskip

$^b$ Faculty of Natural Sciences and Mathematics, University of Maribor, Slovenia\\
\medskip

$^c$ Institute of Mathematics, Physics and Mechanics, Ljubljana, Slovenia\\
\medskip

$^d$ Department of Applied Mathematics, Faculty of Mathematical Sciences,\\
Ferdowsi University of Mashhad, P.O.\ Box 1159, Mashhad 91775, Iran\\
{\tt m$\_$tavakoli@um.ac.ir}

\end{center}

\begin{abstract}
If $S=\{v_1,\ldots, v_k\}$ is an ordered subset of vertices of a connected graph $G$ and $e$ is an edge of $G$, then the vector $r_G(e|S) = (d_G(v_1,e), \ldots, d_G(v_k,e))$ is the edge metric $S$-representation of $e$. If the vertices of $G$ have pairwise different edge metric $S$-representations, then $S$ is an edge metric generator for $G$. The cardinality of a smallest edge metric generator is the edge metric dimension ${\rm edim}(G)$ of $G$. A general sharp upper bound on  the edge metric dimension of hierarchical products $G(U)\sqcap H$ is proved. Exact formula is derived for the case when $|U| = 1$. An integer linear programming model for computing the edge metric dimension is proposed. Several examples are provided which demonstrate how these two methods can be applied to obtain the edge metric dimensions of some applicable graphs.
\end{abstract}


\medskip\noindent
{\bf AMS Subj.\ Class:} 05C12; 05C76

\section{Introduction}
\label{sec:intro}

Graphs considered in this paper are connected, finite, and simple. If $G$ is a graph and $u, v\in V(G)$, then $d_G(u,v)$ denotes the shortest-path distance between $u$ and $v$.  If $S=\{v_1,\ldots,v_k\}$ is an ordered  subset of $V(G)$, then the {\em metric $S$-representation} of a vertex $u\in V(G)$ is the vector $r_G(u|S) = (d_G(v_1,u), \ldots, d_G(v_k,u))$. The set $S$ {\em distinguishes} vertices $u$ and $v$ if $r_G(u|S)\neq r_G(v|S)$ and 
$S$ is a {\em metric generator} for $G$ if each pair of vertices of $G$ is distinguished by $S$. A metric generator of smallest cardinality is called a {\em metric basis} for $G$, its order being the {\em metric dimension} ${\rm dim}(G)$ of $G$.

The sources for the metric dimension are papers~\cite{Harary, Slater}. Afterwards the concept was studied in depth, classical references include~\cite{bailey-2011, caceres-2007, Chartrand-2000}, papers dealing with applications of the metric dimension in modeling of real world problems include~\cite{johnson-1993, khuller-1996}, while for some of the recent developments we refer to~\cite{vetrik-2020, zhang-2020}. Several variations of the concept  were also studied such as the local metric dimension~\cite{okamoto-2010}, independent resolving sets~\cite{Chartrand-2003}, strong resolving sets~\cite{Oellermann-Peters-Fransen}, and $k$-metric generators~\cite{Estrada-Moreno}. Distinguishing edges instead of vertices seems an utmost natural variation, hence it comes as a surprise that the edge metric dimension was introduced only recently in~\cite{kelenc-2018} as follows. 

Let $G$ be a graph. If $u\in V(G)$ and $xy\in E(G)$, then the distance $d_G(u,xy)$ between $u$ and $xy$ is $\min\{d_G(u,x),d_G(u,y)\}$. If $S=\{v_1,\ldots,v_k\}$ is an ordered subset of $V(G)$, then the {\em edge metric $S$-representation} of an edge $e\in E(G)$ is the vector 
$$r_G(e|S) = (d_G(v_1,e), \ldots, d_G(v_k,e))\,.$$
$S$ is an {\em edge metric generator} for $G$ if the edges of $G$ have pairwise different edge metric $S$-representations. A smallest edge metric generator is an {\em edge metric basis} for $G$, its cardinality is the {\em edge metric dimension} $\edim(G)$ of $G$. 

When someone thinks of a smart city, an intelligent transportation system (ITS) may quickly come to mind. Self-driving cars will probably soon play a crucial role in an ITS. Clearly, a self-driving car needs to determine its position on the city's streets uniquely, hence each street needs a code which uniquely determines its location. If we represent the city with a graph $G$, where the edges of $G$ correspond to streets, then an edge metric generator of $G$ provides unique codes for the streets. 

The seminal paper~\cite{kelenc-2018} on the edge metric dimension brings a wealth of results, including a proof that the problem of finding the edge metric dimension of a graph is NP-hard and some approximation results for the invariant. It is also shown that ${\rm dim}(G)$ and $\edim(G)$ are in general incomparable, but  it seems that in most cases ${\rm dim}(G) \le \edim(G)$ holds. In a subsequent paper~\cite{zubrilina-2018} several problems from~\cite{kelenc-2018} are answered, in particular, a classification of the graphs $G$ of order $n$ for which $\edim(G) = n-1$ holds is given. These graphs were also investigated in~\cite{zhu-2019} where a polynomial algorithm is developed for their recognition. Papers~\cite{adawiyah-2019, filipovic-2019, yang-2019,  zhang-2020b} determine the edge metric dimension for some families of graphs. Finally, in~\cite{peterin-2019} the edge metric dimension of the join of graphs, the lexicographic product of graphs, and the corona product of graphs is reported. 

In the next section we study the edge metric dimension of hierarchical products $G(U)\sqcap H$ of graphs. We prove a general sharp upper bound on  $\edim(G(U)\sqcap H)$ and an exact result for the case when $|U| = 1$. Earlier known results on the corona product of graphs can be deduced from these results. In Section~\ref{sec:ILP} we propose an integer linear programming model for computing the edge metric dimension. In the final section several examples are provided that demonstrate how the methods proposed in the previous two sections can be applied to obtain the edge metric dimension of some interesting graphs, notably from mathematical chemistry.

To conclude the introduction we extend (edge) metric generators to vertex and edge subsets as follows. If $X\subseteq V(G)$, then $S\subseteq V(G)$ is a {\em metric generator for $X$} if the vertices from $X$ have pairwise different metric $S$-representations. A smallest metric generator for $X$ is a {\em metric basis} for $X$, its cardinality being the {\em metric dimension} ${\rm dim}_G(X)$ {\em for $X$}. In this notation, ${\rm dim}_G(V(G)) = {\rm dim}(G)$. Similarly, if $F\subseteq E(G)$, then $S\subseteq V(G)$ is an {\em edge metric generator for $F$} if the edges from $F$ have pairwise different edge metric $S$-representations. A smallest edge metric generator for $F$ is an {\em edge metric basis} for $F$, its cardinality is the {\em edge metric dimension} $\edim_G(F)$ {\em for $F$}. So $\edim_G(E(G)) = \edim(G)$.

\section{Hierarchical products}
\label{sec:hierarchical-products}

In this section we consider the edge metric dimension of the hierarchical product of graphs and mention in passing that the metric dimension and the fractional metric dimension of these products were studied in~\cite{feng-2013}, and the local metric dimension in~\cite{klavzar-2020}. 

If $G$ and $H$ are graphs and $U\subseteq V(G)$, then the {\em hierarchical product} $G(U)\sqcap H$ of $G$ and $H$ (with respect to $U$) has the vertex set $V(G)\times V(H)$ and the edge set  
$$\{(g,h)(g',h'):\ gg'\in E(G), h=h'\} \cup \{(g,h)(g',h'):\ g=g'\in U, hh'\in E(H)\}\,.$$
Note that $G(U)\sqcap H$ contains $n(G)$ subgraphs isomorphic to $G$, they are called {\em $G$-layers}. Similarly, $G(U)\sqcap H$ contains $|U|$ subgraphs isomorphic to $H$, these are {\em $H$-layers}. The operation $\sqcap$ (for two and also more factors) was in the seminal paper~\cite{barriere-2009} named the {\em generalized hierarchical product}, here we follow the reasonable suggestion from~\cite{anderson-2017} to simplify the naming to the hierarchical product. 

If $U\subseteq V(G)$ and $u,v\in V(G)$, then we say that a $u,v$-walk $W$ is a {\em $u,v$-walk through $U$} if $W$ is an $u,v$-walk in $G$ that contains some vertex of $U$, where the latter vertex could be one of $u$ and $v$. With $d_{G(U)}(u,v)$ we denote the length of a shortest $u,v$-walk through $U$. With this notation we can state the following fundamental observation from~\cite{barriere-2009}.

\begin{proposition} 
\label{prp:distance}
If $G$ is a graph with $U\subseteq V(G)$ and $H$ is a graph, then 
$$d_{G(U)\sqcap H}((g,h),(g',h'))=\begin{cases}
d_{G(U)}(g,g')+d_H(h,h'); &  h\neq h',\\
d_G(g,g'); & h=h'.
\end{cases}$$
\end{proposition}

To state our results, we need some more preparation. If $v$ is a vertex of a graph $G$ and $k\in {\mathbb N}_0$, then let $E_G(v,k)$ be the set of edges of $G$ that are at distance $k$ from $v$, that is,  
$$E_G(v,k) = \{e\in E(G):\ d_G(v,e) = k\}\,.$$
If $F\subseteq E(G)$ and $|F|\ge 2$, then we say that $X\subseteq V(G)$ is an {\em equidistant discriminator for} $F$, if $X$ is an edge metric generator for $F$. In the case when $|F|\le 1$, we define $\emptyset$ to be the only  equidistant discriminator for $F$. We analogously define equidistant discriminators for vertex subsets of $G$. With this terminology we set    
$$\edim(G(U)) = \min\{ \mid\! \bigcup_{u\in U\atop k\ge 0} S_G(u,k)\! \mid\ :\ S_G(u,k)\ {\rm equidistant\ discriminator\ for}\ E_G(v,k) \}\,.$$
That is, $\edim(G(U))$ is the cardinality of a smallest set of vertices which distinguish all pairs of edges that are equidistant from some vertex from $U$. In addition, we set 
$$\edim^+(G(U)) = \min \left\{ \mid \left(\bigcup_{u\in U, k\ge 0} S_G(u,k)\right) \cup S_G(U) \mid \right\}\,,$$
where the minimum is taken over all equidistant discriminators $S_G(u,k)$ for $E_G(u,k)$ and over all equidistant  discriminators $S_G(U)$ for $U$. After this preparation we can state the following bound. 

\begin{theorem}
\label{t1}
If $G$ and $H$ are graphs and $U\subseteq V(G)$ with $|U|>1$, then
$$\edim(G(U)\sqcap H)\leq n(H) (\edim^+(G(U))+1)\,.$$
\end{theorem}

\proof
Note that the assumption $|U|>1$ implies that also $n(G) > 1$. To simplify the notation, set $X =  G(U)\sqcap H$ for the rest of the proof. Let $S_G(u,k)$, $u\in U$, $k\ge 0$, be equidistant discriminators for $E_G(u,k)$, and $S_G(U)$ be an equidistant discriminator for $U$ which together realize $\edim^+(G(U))$. Set $S^T(G) = \bigcup_{u\in U, k\ge 0} S_G(u,k)$ and let $S=(S^T(G)\cup S_G(U))\times V(H)$. Select further a vertex $w\in U$ and set $S'=\{(w,y):\ y\in V(H)\}$.
We claim that $S\cup S'$ is an edge metric generator for $X$. For this sake let $e$ and $f$ be arbitrary, different edges of $X$, and consider the following cases. 

\medskip\noindent
{\bf Case 1}: $e$ and $f$ are both in $H$-layers.\\
Suppose first that $e$ and $f$ are in the same $H$-layer, say $e = (g,h)(g,h')$ and $e = (g,h'')(g,h''')$. If $\{h,h'\}\cap \{h'',h'''\} = \emptyset$, that is, if $e$ and $f$ are not adjacent, then consider an arbitrary vertex $v\in S^T(G)$ and note that the vertex $(v,h)\in S$ distinguishes $e$ and $f$. In the second subcase  suppose that $e$ and $f$ are adjacent, say $h = h''$. Now the vertex $(v,h')\in S$ distinguishes $e$ and $f$. 

Suppose second that $e$ and $f$ are in different $H$-layers, say $e = (g,h)(g,h')$ and $e = (g',h'')(g',h''')$, where $g\ne g'$. Select a vertex $v\in S_G(U)$ which distinguishes $g$ and $g'$, that is, $d_G(v,g)\neq d_G(v,g')$. In the first subcase suppose that $\{h,h'\}\cap \{h'',h'''\} \ne \emptyset$, say $h=h''$. It this case the vertex $(v,h)\in S$ distinguishes $e$ and $f$. Suppose next that $\{h,h'\}\cap \{h'',h'''\} =  \emptyset$. Then  $d_X((v,h),e)=d_X((v,h), (g,h))$ and we may without loss of generality assume that $d_X((v,h),f)=d_X((v,h), (g',h''))$. If $d_X((v,h),e) \ne d_X((v,h),f)$, then $(v,h)$ distinguishes $e$ and $f$. Suppose next that $d_X((v,h),e) = d_X((v,h),f)$. Then $d_{G(U)}(v,g)=d_{G(U)}(v,g')+d_H(h,h'')$. Thus $d_{G(U)}(v,g)+d_H(h,h'')>d_{G(U)}(v,g)=d_{G(U)}(v,g')+d_H(h,h'')>d_G(v,g')+d_H(h'',h'')$ and so $d_X((v,h''),e) >  d_X((v,h''),f)$. Therefore, the vertex $(v,h'')\in S$ distinguishes $e$ and $f$.

\medskip\noindent
{\bf Case 2}: $e$ and $f$ are  are both in $G$-layers.\\
Let $e=(g,h)(g',h)$ and $f=(g'',h')(g''',h')$. 
First, we check the case $h=h'$. In this case, if there exists a vertex $v\in U$ such that $d(v,gg')=d(v,g''g''')$, then there exists a vertex $u\in S^T(G)$ for which $d_G(u,gg')\neq d_G(u,g''g''')$ holds. Consequently, for the vertex $(u,h)\in S$ we have $d_X((u,h),e) \neq d_X((u,h),f)$. 
Otherwise, the vertex $(w,h)\in S'$ distinguishes $e$ and $f$.\\
Now we investigate the case $h\neq h'$. In this case, if $d_X((u,h),e)= d_X((u,h),f)$ for each $(u,h)\in S$, then, again by a similar argument as applied in Case 1, there exists $(u,h')$ in $S$ that $d_X((u,h'),e)\neq d_X((u,h'),f)$.

\medskip\noindent
{\bf Case 3}: $e$ is in a $G$-layer and $f$ is in a $H$-layer.\\ 
Let $e=(g,h)(g',h)$ and $f=(g'',h')(g'',h'')$.
The case $h\notin \{h',h''\}$ can be proved by a similar technique as used in Case 1 and so we will only check the case $h=h'$.
If $d_X((u,h),e)= d_X((u,h),f)$ for each $(u,h)\in S$, then $\min\{d_G(g,u),d_G(g',u)\}=d_G(u,g'')$ for each $u\in S^T(G)\cup S_G(U)$.
Therefore, $d_X((u,h''),f)=d_X((u,h''), (g'',h''))=d_G(u,g'')<d_X((u,h''),e)$. So we have detected the vertex $(u,h'')\in S$ such that $d_X((u,h''),e) \neq d_X((u,h''),f)$.    

We conclude that every pair of edges from $X$ is distinguished by a vertex of $S\cup S'$,
and consequently  $\edim(X) \leq n(H) (\edim^+(G(U))+1)$.
\qed

Consider the cases in which $|U|=1$ or $U\cap (S^T(G)\cup S_G(U))\neq\emptyset$, where $S^T(G)$ and  $S_G(U)$ are defined as in the proof of Theorem~\ref{t1}. If there exists $z\in U\cap (S^T(G)\cup S_G(U))$, then replacing the vertex set $S'$ in the proof of Theorem~\ref{t1} with the set$\{(z,y):\ y\in V(H)\}$  we can prove along with the lines of the proof that $S$ is an edge metric generator for $G(U)\sqcap H$. In the case when $|U|=1$, Cases 1 and~3 from the proof are valid for $G(U)\sqcap H$, and we do not need the vertices of $S'$ in Case 2. Summarizing this discussion we have the following fact. If $|U|=1$ or $U\cap (S^T(G)\cup S_G(U))\neq\emptyset$, then 
\begin{equation}
\label{eq:1}
\edim(G(U)\sqcap H)  \leq n(H)\,\edim^+(G(U))\,. 
\end{equation}

\begin{figure}[ht!]
\centerline{ \includegraphics[scale=0.8]{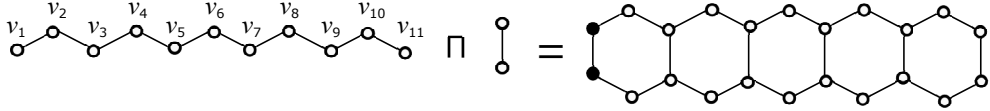}}
\caption{$P_{11}(U)\sqcap P_2$, where $U=\{v_1,v_3,v_5,v_7,v_9,v_{11}\}$.}
\label{fig3}
\end{figure}

Consider $P_{11}(U)\sqcap P_2$, where $V(P_{11}) = \{v_1,\ldots, v_{11})$ and $U=\{v_{2k-1}:\ k\in [6]\}$, see Fig.~\ref{fig3}. 

Then we can select $S_{P_{11}}(U) = \{v_1\}$ and $S^T(P_{11})=\{v_1\}$, and so $\edim^+(P_{11}(U))=1$. From~\eqref{eq:1} we then infer that $S=(S_{P_{11}}(U)\cup S^T(G))\times V_H=\{v_1\}\times V_H$ (black vertices in the figure) is an edge metric generator for $P_{11}(U) \sqcap P_2$, hence ${\rm edim}(P_{11}(U) \sqcap P_2) \leq 2$. From~\cite[Remark 1]{kelenc-2018} we know that  ${\rm edim}(G)=1$ if and only if $G$ is a path, therefore we conclude that ${\rm edim}(P_{11}(U) \sqcap P_2)= 2$. This demonstrates that~\eqref{eq:1} is sharp.

We now focus on hierarchical products $G(U)\sqcap H$, where $|U| = 1$. If $U = \{u\}$, then we simplify the notation $G(\{u\})$ to $G(u)$. If $G$ is a path and $u$ its end vertex, then we say that $G(u)$ is a {\em rooted path}. 

\begin{theorem}
\label{t2}
If $X = G(u)\sqcap H$, where $G(u)$ is not a rooted path, and $n(H)\ge 2$, then 
$${\rm edim}(X) = n(H) \cdot \edim(G(u))\,.$$
\end{theorem}

\proof
From~\eqref{eq:1} we know that ${\rm edim}(X) \le n(H) \cdot \edim^+(G(u))$. Since $|U|=1$, the  equidistant discriminator for $U$ is the empty set, that is, $S_G(U) = \emptyset$, and consequently $ \edim(X)\leq n(H)\cdot \edim(G(u))$. 

Let $S^T(G) = \bigcup_{u\in U, k\ge 0} S_G(u,k)$, where $S_G(u,k)$ are equidistant discriminators for $E_G(u,k)$  that realize $\edim(G(U))$.  Then we know that $S=S^T(G)\times V(H)$ is an edge metric generator for $X$. We wish to show that $|S| = \edim(X)$ and assume  by way of contradiction that there is an edge metric generator $S'$ for $X$ such that $|S'|<|S|$. By the pigeonhole principle there exists a $G$-layer of $X$, denote it with $G_h$ (here $h$ is the vertex of $H$ to which the $G$-layer corresponds), such that $|S'\cap V(G_h)| < |S\cap V(G_h)|$. Let $S'_h = S' \cap V(G_h)$ and note that $S'_h$ is not an equidistant discriminator for $E(G_h)$ and $|S'_h|<|S^T(G)|$. Hence there exist $k\ge 0$ and edges $e, f \in E(G_h)\cap \{(g,h)(g',h):\ gg'\in E_G(u,k)\}$ such that $d_X(x,e) = d_X(x,f)$ holds for each vertex $x \in S'_h$. Since $d_X((u,h),e)=d_{X}((u,h),f)$, it follows that the equality $d_X(v,e) = d_X(v,f)$ holds also for each $v\in S'\setminus V(G_h)$. But this means that $S'$ is not an edge metric generator for $X$, a contradiction. 
\qed

Note that the rooted paths $G(u)$ are the only graphs for which $S^T(G) = \emptyset$ and consequently $\edim(G(u)) = 0$. This is the reason that the rooted paths are excluded in Theorem~\ref{t2}. 

To conclude the section we consider the corona product of graphs. Recall that the {\em corona product} $G\odot H$ of graphs $G$ and $H$ is obtained from the disjoint union of $G$ and $n(G)$ copies of $H$, by joining by an edge every vertex from the $i^{\rm th}$ copy of $H$ with the $i^{\rm th}$ vertex of $G$. (See~\cite{Corona} for more information on this product.) The key observation is that  
$$G\odot H = (H+v)(v) \sqcap G\,,$$
where $H+v$ denotes the join of $H$ and the one vertex graph with the vertex $v$. More precisely, $H+v = H+K_1$, where $V(K_1) = \{v\}$.  Then Theorem~\ref{t2} implies that 
$$\edim(G\odot H) = \edim((H+v)(v) \sqcap G) = n(G)\cdot \edim((H+v)(v))\,.$$
From here it is not difficult to deduce~\cite[Theorem 4.1]{zubrilina-2018} which determines the edge metric dimension for the join of $K_1$ and an arbitrary graph, and~\cite[Theorem 6]{peterin-2019} that determines the edge metric dimension of corona products of nontrivial graphs. We can reformulate and combine these two results as follows. Let $\mathcal{F}$ be the family of graphs consisting of all graphs $G$ such that that $d_G(v,e)\leq 1$ holds for each $e\in E(G)$ and each $v\in V(G)$. 

\begin{theorem}
\label{t3}
If $G$ is a connected graph, and $H$ is a graph with more than one vertex, then
$${\rm edim}(G\odot H)=
\begin{cases}
n(H); &  n(G) = 1\ \text{and}\ H\in \mathcal{F}, \\
n(G)(n(H)-1); & \text{otherwise}.
\end{cases}$$
\end{theorem}

\section{Integer linear programming model}
\label{sec:ILP}

An integer linear programming model, ILPM for short, for finding the metric dimension and a metric basis for a graph has been presented in~\cite{Chartrand-2000}. Following this approach we introduce an ILPM for finding the edge metric basis for a given graph as follows. 

Let $G$ be a graph with $V(G) = \{v_1,\ldots,v_n\}$ and $E(G) = \{e_1,\ldots,e_m\}$.  Let $D_G=[d_{ij}]$ be an $m\times n$ matrix, where $d_{ij}=d_G(e_i,v_j)$ for $i\in [m]$ and $j\in [n]$. For $x_i \in \{0,1\}$, $i\in [n]$, define the function $F$ by 
\[F(x_1,\ldots, x_n)=x_1+\cdots+x_n\,, \]
and minimize $F$ subject to the constraints
\[|d_{i1}-d_{j1}|x_1+|d_{i2}-d_{j2}|x_2 + \cdots + |d_{in}-d_{jn}|x_n>0,\ 1 \le i < j \le m\,.\]
Then note that if $x'_1, \ldots , x'_n$ is a set of values for which $F$ attains its minimum, then $W = \{v_i:\ x'_i =1\}$ is an edge metric basis for $G$. 

For example, consider $K_3$ with the vertex set $\{v_1,v_2,v_3\}$ and edges $e_1=v_1v_2$, $e_2=v_2v_3$, and $e_3=v_1v_3$. Then $D_{K_3}=\begin{pmatrix}
0 & 0& 1\\
1 & 0 & 0\\
0 & 1 & 0
\end{pmatrix}$. Thus, minimize $F(x_1,x_2,x_3)=x_1+x_2+x_3$ subject to the  constraints $x_1+x_3 > 0$, $x_2 + x_3 > 0$, $x_1+x_2 > 0$, $x_1, x_2, x_3 \in \{0, 1\}$. Then $F$ attains its minimum for $x_1=0$, $x_2=1$, and $x_3=1$, hence $W=\{v_2,v_3\}$ is an edge metric basis for $K_3$.

\section{Applications}
\label{sec:applications}

In this section we demonstrate how the results from the previous sections can be applied to compute the edge metric dimension of interesting graphs.

\begin{figure}[ht!]
\centerline{ \includegraphics[scale=0.25]{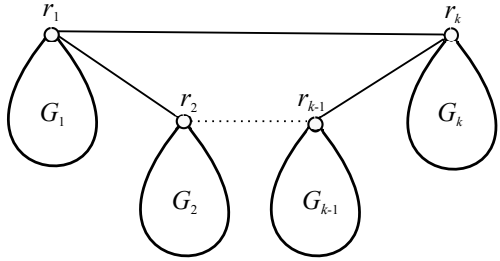}}
\caption{The bridge-cycle graph $BC(G_1,\ldots,G_k; r_1,\ldots,r_k)$.}
\label{fig7}
\end{figure}

Let $G_1,\ldots,G_k$ be rooted graphs with respective root vertices $r_1,\ldots,r_k$. The {\em bridge-cycle graph} $BC(G_1,\ldots,G_k; r_1,\ldots,r_k)$ is the graph obtained from the disjoint union of $G_1,\ldots,G_k$ by joining the vertices $r_i$ and $r_{i+1}$ for all $i\in [r-1]$ and adding the edge $r_1r_k$,  see Fig.~\ref{fig7}. 

If $G_1=\cdots=G_k=G$ and $r=r_1$, where $G(r)$ is not a rooted path, then we infer that $BC(G_1,\ldots,G_k; r_1,\ldots,r_k)\cong G(r)\sqcap C_k$. Theorem \ref{t2} then implies that
 $${\rm edim}(BC(G_1,\ldots,G_k; r_1,\ldots,r_k))={\rm edim}(G(r)\sqcap C_k)=k\cdot {\rm edim}(G(r))\,.$$

\begin{figure}[ht!]
\centerline{ \includegraphics[scale=0.45]{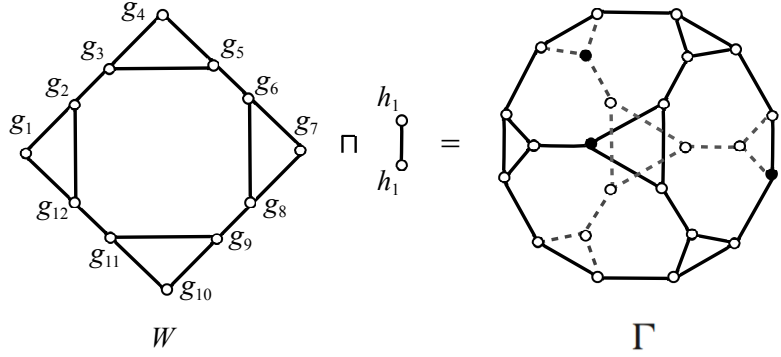}}
\caption{$W(U)\sqcap P_2=\Gamma$ where $U=\{g_1,g_4,g_7,g_{10}\}$.}
\label{fig4} 
\end{figure}

The examples from the rest of this section come from chemical graph theory. Consider first the (molecular) graph of truncated cube, it is denoted by $\Gamma$ and drawn in Fig.~\ref{fig4}. 

As the figure shows, $\Gamma$ is isomorphic to the hierarchical product $W(U)\sqcap P_2$ (see the figure for $W$),  where $U=\{g_1,g_4,g_7,g_{10}\}$. Then, by the proof of Theorem~\ref{t1}, $S=\{g_1,g_6\}\times V(H) = \{(g_1,h_1),(g_1,h_2),(g_6,h_1),(g_6,h_2)\}$ because $S^T(W) = S_W(U) = \{g_1,g_6\}$.
Thus ${\rm edim}^+(W(U))=2$, and so by~\eqref{eq:1}, ${\rm edim}(\Gamma)={\rm edim}(W(U)\sqcap P_2)\leq 4$. On the other hand, the exact value of ${\rm edim}(\Gamma)$ computed by the ILPM from Section~\ref{sec:ILP} is equal to $3$. The black vertices from Fig.~\ref{fig4} form an edge metric bases of $\Gamma$ found by the ILPM.

\begin{figure}[ht!]
\centerline{ \includegraphics[scale=0.35]{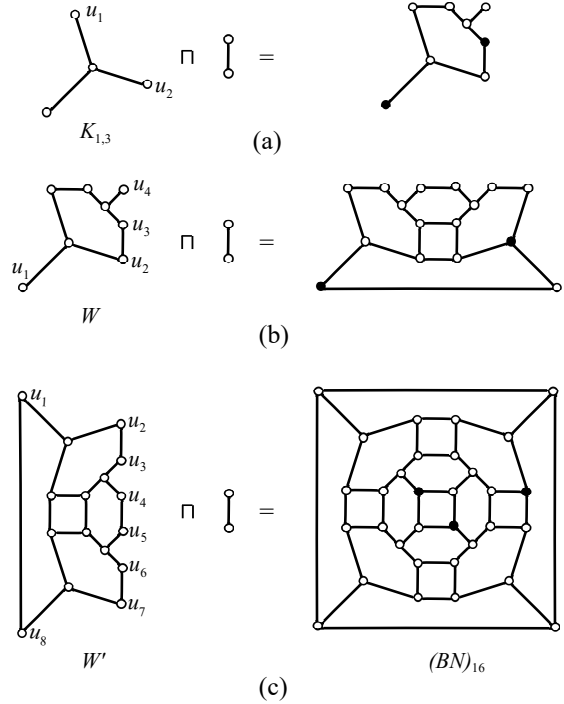}}
\caption{ (a) $K_{1,3}(U)\sqcap P_2$ where $\{u_1,u_2\}$\ (b) $W(U)\sqcap P_2$ where $\{u_1,\ldots,u_4\}$\ (c) $W'(U)\sqcap P_2$ where $\{u_1,\ldots,u_8\}$.}
\label{fig2}
\end{figure}

Continuing with examples from chemical graph theory, recall that a {\em fullerene} is a plane, $3$-connected, cubic graph with only pentagonal and hexagonal faces. The literature on fullerenes is huge, see for instance~\cite{fullerene1} for more informations about their electronic and structural properties and the recent survey~\cite{andova-2016}. More generally, the term fullerene is also used for such graphs where other lengths of faces are present, cf.~\cite{shi-2019, zhao-2018}. For instance, the graph $(BN)_{16}$ from Fig.~\ref{fig2} is an example of a $(4,6)$-fullerene. 

By Fig.~\ref{fig2} and Theorem \ref{t1}, we have ${\rm edim}(K_{1,3}(U)\sqcap P_2)\leq 4$, ${\rm edim}(W(U)\sqcap P_2)\leq 4$, and ${\rm edim}((BN)_{16})={\rm edim}(W'(U)\sqcap P_2)\leq 4$. On the other hand, using the ILPM from the previous section we get ${\rm edim}(K_{1,3}(U)\sqcap P_2)= 2$, ${\rm edim}(W(U)\sqcap P_2)= 2$, and ${\rm edim}((BN)_{16})= 3$. The black vertices show form the edge metric bases found by the ILPM.

\end{document}